\documentclass[12pt]{article}\usepackage{amsmath,epsfig,amssymb,amsfonts,amsthm,epsfig,verbatim,bm,hyperref}
\usepackage{color,graphicx,lscape,longtable}
\usepackage[authoryear]{natbib}
\usepackage{tikz}

\newcommand{\indep}{\;\, \rule[0em]{.03em}{.67em} \hspace{-.25em}
\rule[0em]{.65em}{.03em} \hspace{-.25em}\rule[0em]{.03em}{.67em}\;\,}

\newtheorem{Th}{\underline{\bf Theorem}}

\newtheorem{Rem}{\underline{\bf Remark}}

\newtheorem{Pro}{\underline{\bf Proposition}}

\newtheorem{Ex}{Example}

\def\rkhsg{\mathcal{H}_{\Gamma}}

\def\bse{\begin{eqnarray*}}
\def\ese{\end{eqnarray*}}
\def\be{\begin{eqnarray}}
\def\ee{\end{eqnarray}}
\def\bsq{\begin{equation*}}
\def\esq{\end{equation*}}
\def\bq{\begin{equation}}
\def\eq{\end{equation}}

\def\cov{{\rm cov}}

\def\real{\mathbb{R}}

\newcommand{\argmin}{\mathop{\rm argmin}}
\newcommand{\argmax}{\mathop{\rm argmax}}

\def\trans{^{\rm T}}

\def\bb{\boldsymbol\beta}

\def\0{{\bf 0}}

\def\v{{\bf v}}

\def\X{{\bf X}}

\def\Z{{\bf Z}}

\def\bq{\begin{equation}}
\def\eq{\end{equation}}

\def\trans{^{\rm T}}

\def\squarebox#1{\hbox to #1{\hfill\vbox to #1{\vfill}}}

\def\bse{\begin{eqnarray*}}
\def\ese{\end{eqnarray*}}
\def\be{\begin{eqnarray}}
\def\ee{\end{eqnarray}}
\def\bsq{\begin{equation*}}
\def\esq{\end{equation*}}
\def\bq{\begin{equation}}
\def\eq{\end{equation}}

\def\trans{^{\rm T}}

\def\boxit#1{\vbox{\hrule\hbox{\vrule\kern6pt\vbox{\kern6pt#1\kern6pt}\kern6pt\vrule}\hrule}}

\begin{document}
\thispagestyle{empty} \baselineskip 22pt
\renewcommand {\thepage}{}
\include{titre}
\pagenumbering{arabic}
\begin{center}
{\Large\bf Functional Inverse Regression in an Enlarged Dimension
Reduction Space}
\end{center}

\baselineskip 17pt

\vspace{5mm}
\begin{center}
{Ting-Li Chen$^1$, Su-Yun Huang$^1$, Yanyuan Ma$^2$ and I-Ping Tu$^1$}\\
\vspace{5mm}

$^1${\it Institute of Statistical Science\\
Academia Sinica\\
Taipei 11529, Taiwan \\
{\rm email: tlchen, syhuang, iping@stat.sinica.edu.tw}}\\

\vspace{5mm}
$^2${\it Department of Statistics\\
The University of South Carolina\\
Columbia, SC 29208\\
{\rm email: yanyuanma@stat.sc.edu}}\\
\end{center}


\begin{quotation}
\noindent {\it Abstract:}
We consider an enlarged dimension reduction space in functional
inverse regression. Our operator and functional analysis based
approach facilitates a compact and rigorous formulation of the
functional inverse regression problem. It also enables us to expand
the possible space where the dimension reduction functions belong.
Our formulation provides a unified framework so that the classical
notions, such as covariance standardization, Mahalanobis distance,
SIR and linear discriminant analysis, can be naturally and smoothly
carried out in our enlarged space. This enlarged dimension reduction
space also links to the linear discriminant space of Gaussian
measures on a separable Hilbert space. \par

\vspace{9pt}
\noindent {\it Key words and phrases:}
Functional inverse regression, functional dimension reduction, functional linearity condition, sliced inverse regression.
\par
\end{quotation}

\section{Introduction}

Traditionally, sufficient dimension reduction problems refer to the
estimation of the space spanned by the columns of $\bb$, where $\bb$
satisfies $Y\indep \X\mid\bb\trans\X$. Here, $\X$ is a
$p$-dimensional covariate vector, which we assume to satisfy
$E(\X)=\0$ for simplicity,  $\bb$ is a $p\times d$ matrix and $Y$ is
a univariate response variable. An equivalent form is
$Y=f(\bb\trans\X, \epsilon)$, where $\epsilon$ is a mean zero random
variable independent of $\X$. By far the most well known procedure
of estimation in this problem is sliced inverse regression (SIR,
 \cite{li1991}), where solving the leading $d$ eigenvectors
of the generalized eigenvalue problem $\Gamma_e\v=\lambda\Gamma\v$
is all one needs to do to obtain the column space of $\bb$. Here
$\Gamma=\cov(\X)$ and $\Gamma_e=\cov\{E(\X\mid Y)\}$. SIR is
constructed under a linearity condition which requires
$E(\X\mid\bb\trans\X)=\Gamma\bb(\bb\trans\Gamma\bb)^{-1}\bb\trans\X$
and is then further developed into a whole class of inverse
regression based methods for dimension reduction. To understand the
inverse regression based methods from a different angle, we can
normalize the covariates through viewing $\Z \equiv\Gamma^{-1/2}\X$
as new covariates and $\bm{\eta}\equiv\Gamma^{1/2}\bb$ as new
dimension reduction matrix. Considering the dimension reduction
problem in terms of $(\Z, Y, \bm{\eta})$ instead of $(\X,Y,\bb)$
enables much simplification and permits clearer exhibition of the
critical operations~\cite{li1991, mazhu2012}.

Dimension reduction problems have been extended from the traditional
regression domain to the functional data analysis domain. See
\cite{jiangetal2014} and references therein. The model considered in
the functional dimension reduction framework is \bse
Y=f(\langle\beta_1,X\rangle_{L_2},\dots,\langle\beta_d,X\rangle_{L_2},\epsilon),
\ese where $Y$ is still a univariate response variable, $X$ is now a
covariate function, $\beta_1, \dots,\beta_d$ are parameter functions
in $L_2(I)$, and $\langle\,\,\,, \,\,\,\rangle_{L_2}$ denotes the
inner product of two functions in the $L_2(I)$ space. Since the
matrix vector product $\bb\trans\X$ in the traditional case can be
expressed as $(\bb_1\trans\X, \dots,\bb_d\trans\X)\trans$ which can
also be viewed as a vector of inner products between the vector
$\bb_k$ and the covariate $\X$, one might think that the extension
to the functional data framework is straightforward. However, there
are many subtleties when finite dimensional quantities are extended
to infinite dimensional ones, such as $\bb_1, \dots, \bb_d$ and $\X$
to $\beta_1(\cdot), \dots, \beta_d(\cdot)$ and $X(\cdot)$. Some
properties we take for granted in finite dimension may not hold
automatically, e.g., some vector norm or the inner product between
vectors may not be finite. If we want to perform the similar
standardization as in the finite dimensional case by forming
$Z(\cdot)\equiv\Gamma^{-1/2} X(\cdot)$ as new covariate function,
and $\eta_j(\cdot)\equiv\Gamma^{1/2}\beta_j(\cdot)$ as new dimension
reduction functions for the functional correspondence of the
variance-covariance matrix (operator) $\Gamma$, not only do we need
to consider the extensions from vectors to functions and matrices to
operators, but also to define a proper normed spaces and their
corresponding requirements on these functions. A careful and
rigorous consideration of these issues will enable less restrictive
models and more flexible estimation. In fact, one of the main
messages of this article is to point out that the requirement of the
parameter functions $\beta_j(\cdot)$'s being $L_2(I)$ in
\cite{jiangetal2014} is too strong and can be relaxed to include
more interesting examples.

During the process of our investigation, we also realize that it is
crucial to formulate the functional dimension reduction problem
properly in order to facilitate the subsequent application of the
existing mathematical tools from functional analysis involving
Reproducing Kernel Hilbert Space (RKHS) and operator theories. To
better prepare for such a task, we summarize some preliminary
results in Section~\ref{sec:pre} and provide an outline of either a
proof or an understanding for each result. In Section~\ref{sec:ex},
we give a few motivating examples, wherein the dimension reduction
functions fall out of the  space required in \cite{jiangetal2014}
and hence cannot be solved under their model. We then present an
extension of the functional dimension reduction model in
Section~\ref{sec:main}, together with some main results. Our
extension works on an enlarged space, so that the classical notion
of SIR in standardized scale can be carried out.

\section{Preliminary}\label{sec:pre}

\subsection{Covariance operators and integral operators}

Without loss of generality, we restrict our attention to functions
defined on $I\equiv[0,1]$. Let the Hilbert space $L_2(I)$ be the
space of functions defined on $I$ and equipped with inner product
given by
\[\langle u,v\rangle_{L_2}=\int_0^1 u(t)v(t)\, dt,~~~ u, v\in
L_2(I).\] Let $\Gamma(s,t)$ be a  continuous bivariate function on $I\times I$.
Then $\Gamma(s,t)$ induces a  linear integral operator, still
written as $\Gamma(s,t)$, where its operation on a function
$u(\cdot)\in L_2(I)$ is defined as
\begin{equation}\label{eq:Gamma_op}
(\Gamma u)(s)=\int_0^1\Gamma(s,t)u(t)dt
 = \langle \Gamma(s,\cdot), u(\cdot)\rangle_{L_2}~~ \mbox{for $u\in L_2(I)$}.
\end{equation}
When $\langle u,\Gamma v\rangle_{L_2} =\langle \Gamma u,
v\rangle_{L_2}$ for all $u, v\in L_2(I)$,  $\Gamma(s,t)$ is said to
be a symmetric linear integral operator. Note that \bse \langle u,
\Gamma v\rangle_{L_2} =\int_0^1 \int_0^1
u(s)\Gamma(s,t)v(t)\,dt\, ds,\\
\langle \Gamma u, v\rangle_{L_2} =\int_0^1 \int_0^1 u(s)\Gamma(t,s)
v(t) \, dt \,ds. \ese Hence, as long as $\Gamma(s,t)$ is symmetric
as a function of $(s,t)$ defined on $I\times I$,  its induced
operator $\Gamma(s,t)$ is also a symmetric operator. When $\langle
u, \Gamma u\rangle_{L_2} \ge 0$ for all $u\in L_2(I)$, $\Gamma$ is
said to be positive semi-definite (or non-negative definite). When
the equality holds if and only if $u=0$ a.s., $\Gamma$ is said to be
positive definite (or strictly positive definite). A positive
(semi-) definite linear integral operator is also known as a
covariance operator. Let ${\cal B}$ denote the unit ball in
$L_2(I)$, i.e., ${\cal B}\equiv\{f\in L_2(I): \|f\|_{L_2}\le 1\}$.
An operator $\Gamma$ defined on $L_2(I)$ that maps to $L_2(I)$ is
said to be compact if the image of the unit ball, $\Gamma({\cal
B})$, is a compact set in $L_2(I)$.

Let $X(t)$, $t\in I$, be a random process with finite second moments
and $Y$ be a univariate random variable. We now consider three
specific bivariate functions and their induced operators, \bse
 \Gamma(s,t)\equiv \cov\{X(s), X(t)\},~~
\Gamma_w(s,t) \equiv E\left[\cov\{X(s),X(t)|Y\}\right], \ese
 and
\bse
 \Gamma_e(s,t)\equiv\cov [E\{X(s)\mid Y\}, E\{X(t)\mid Y\}].
\ese It is easy to verify that $\Gamma(s,t)$, $\Gamma_w(s,t)$ and
$\Gamma_e(s,t)$ are all symmetric bivariate functions and
$\Gamma(s,t)=\Gamma_w(s,t)+\Gamma_e(s,t)$. We further assume
$\Gamma(s,t)$, $\Gamma_w(s,t)$, $\Gamma_e(s,t)$ to be continuous.
The continuity of  functions $\Gamma(s,t)$, $\Gamma_w(s,t)$ and
$\Gamma_e(s,t)$ on $I\times I$ implies they are square integrable,
and hence the continuity guarantees that, $\Gamma(s,t)$,
$\Gamma_w(s,t)$ and $\Gamma_e(s,t)$ are compact operators on
$L_2(I)$
  \cite{lax2002} (Chapter 22, Theorem 4).
The definitions of $\Gamma(s,t)$, $\Gamma_w(s,t)$ and
$\Gamma_e(s,t)$ also ensure that, they are positive semi-definite.
Mercer's Theorem \cite{lax2002} (Chapter~30, Theorem~11) then
implies that they have discrete spectra. Taking $\Gamma(s,t)$ for
instance, it can be expanded in a uniformly convergent series of
eigenvalues and eigenfunctions \be\label{eq:gamma}
\Gamma(s,t)=\sum_{i=1}^q\xi_i \phi_i(s)\phi_i(t), ~~ q\le\infty,\ee
which we sometimes write in short as
\[\Gamma = \sum_{i=1}^q \xi_i \phi_i\otimes\phi_i\trans .\]
Here $\xi_1\ge \xi_2\ge\dots\ge\xi_q>0$ are decreasing positive
values. If $\Gamma(s,t)$ is strictly positive definite, then
$q=\infty$ and $\{\phi_i(\cdot)\}_{i=1}^q$ form a complete
orthonormal basis for $L_2(I)$. The above result critically relies
on the strictly positive definiteness of $\Gamma$. Without the
assumption of $\Gamma$ being strictly positive definite, we can
still decompose $\Gamma(s,t)$ as in (\ref{eq:gamma}), and the
corresponding $\{\phi_i(\cdot)\}_{i=1}^q$, $q \le \infty$,  always
form a complete orthonormal basis for $R(\Gamma)$, the range of
$\Gamma$.  However, $R(\Gamma)\subsetneq L_2(I)$, when $\Gamma$ is
not strictly positive definite. We further outline the following
results which are relevant to the functional inverse regression
study.
\begin{Pro}\label{prop:1}
A continuous, symmetric, positive (semi-) definite integral operator
$\Gamma(s,t)=\sum_{i=1}^q \xi_i\phi_i(s)\phi_i(t)$ is a trace-class
operator, i.e., \bse
 \sum_{i=1}^q \xi_i <\infty.
\ese
\end{Pro}
\begin{proof}
 Because $\Gamma(s,t)$ is a continuous function on
  $I\times I$, for $s=t$, $f(t)\equiv \Gamma(t,t)$ is a continuous
  function of $t$ in $I$, thus is integrable. Hence,
$\sum_{i=1}^q\xi_i = \int f(t)dt <\infty$.
\end{proof}

\begin{Pro}\label{pro:2}
For any positive (semi-) definite operator $\Gamma(s,t)=\sum_{i=1}^q
\xi_i\phi_i(s)\phi_i(t)$, there exists a mean zero random process
$X(s)$ satisfying $\int E\{X^2(s)\}ds<\infty$ such that
$\Gamma(s,t)=\cov\{X(s),X(t)\}$ and \bse
 X(s) =\sum_{i=1}^q A_i\phi_i(s),
\ese
 where $A_i$'s are independent random variables with mean zero and
variances $\xi_i$'s.
\end{Pro}
\begin{proof} For $i=1,2,\dots, q$, let $A_i=\xi_i^{1/2}Z_i$, where $Z_i$'s
are independent standard normal random variables. Obviously the
resulting $X(s)$ is a mean zero process that satisfies
$\cov\{X(s),X(t)\}=\Gamma(s,t)$. In addition, $\int E\{X^2(s)\}ds
=\sum_{i=1}^q\xi_i<\infty$.
\end{proof}

Note that, in our construction of the Gaussian process in the above
proof, the sample path $X(\cdot|\omega)$ may not be in $L_2(I)$ for
a given realization $\omega$. However, $\int E\{X^2(s)\}ds<\infty$
ensures that the probability of this kinds of $\omega$ is 0. That
is, $X(\cdot|\omega) \in L_2(I)$ almost surely. In the following, we
may simply use $X\in L_2(I)$ to denote that $X(\cdot|\omega) \in
L_2(I)$ almost surely.

\subsection{RKHS relevant for functional inverse regression}

Let $\rkhsg$ be the RKHS generated by $\Gamma(s,t)$. Specifically,
\bse \rkhsg \equiv {\rm closure}{\Big\{\sum_{i=1}^q
\Gamma(s,t_i)\alpha_i: q\in{\mathbb N},\alpha_i\in{\mathbb R},
t_i\in[0,1]\Big\}}, \ese where the closure is taken with respect to
the norm induced by the following inner product \bse \langle
\Gamma(s,\cdot),\Gamma(t,\cdot)\rangle_{\rkhsg}=\Gamma(s,t). \ese
Note that $\rkhsg$ is a {\it proper subset} of $L_2(I)$. For $f\in
\rkhsg\subset L_2(I)$, $f$ has the expansion \bse f(t)=\sum_i
f_{i}\phi_{i}(t),~~\mbox{where $f_{i}=\langle
f,\phi_{i}\rangle_{L_2}$.} \ese
 In addition to its $L_2$-norm
defined as $\|f\|_{L_2}=\sum_i f_i^2$,  the $\rkhsg$-norm is given
by \bse \Vert f\Vert_{\rkhsg}^2 = \sum_i\frac{f_{i}^2}{\xi_{i}}.
\ese For $u,v\in\rkhsg$, the $\rkhsg$-inner product is given by
\begin{equation}\label{rkhs_norm}
 \langle u,v\rangle_{\rkhsg}=\sum_i \frac{u_iv_i}{\xi_i},
\end{equation}
where $u(t)=\sum_i u_i\phi_i(t)$, $v(t)=\sum_i v_i\phi_i(t)$.

\section{Motivating Examples}\label{sec:ex}

Throughout our development of a rigorous framework for functional
inverse regression, we set up a space,
\[R(\Gamma^{-1/2})\equiv \left \{f: f=\sum_{i=1}^\infty
f_i\phi_i, f_i\in\real~\mbox{such that~}
\sum_i\xi_if_i^2<\infty\right\} \supsetneq L_2(I), \] which is the
range space of the operator $\Gamma^{-1/2}$ and is larger than
$L_2(I)$. Below we give a few examples, wherein the dimension
reduction functions fall out of $L_2(I)$ and reside in
$R(\Gamma^{-1/2})$. These examples motivate us to consider an
enlarged space for functional dimension reduction. Interestingly,
this enlarged space $R(\Gamma^{-1/2})$ is the space considered by
Grenander~\cite{grenander1950} and Rao and
Varadarajan~\cite{rao1963} in the study of linear discriminant
analysis of Gaussian measures on a separable Hilbert space.

\begin{Ex}[Binary response]{\rm
Let $Y$ be a binary random variable having probabilities
$P(Y=1)=P(Y=-1)=\frac12$, and let $\{\psi_i\}_{i=1}^\infty$ be a
complete orthonormal basis for $L_2(I)$. Given $Y=y$, consider \bse
X_y (t) = \alpha y\sum_{i=1}^\infty \frac 1{i^{2+\delta}}  \psi_i(t)
 + \sum_{i=1}^\infty \frac 1i \, Z_i\, \psi_i(t),~~ t\in I,
\ese where $0<\delta\le 1/2$ and $\alpha$ is some scalar that
controls the separation of two groups. Here $Z_i$'s are independent
standard normal random variables that are independent of $Y$. Let
$\Gamma_e$ be the between-group covariance and $\Gamma_w$ be the
within-group covariance. Then $\Gamma=\Gamma_e+\Gamma_w$. We can
easily calculate the within-group covariance function as
\[
\Gamma_w(s,t)\equiv E\left[\cov\{X(s), X(t)|Y\}\right]
=\sum_{i=1}^\infty \frac 1 {i^2} \psi_i(s) \psi_i(t),\] and the
between-group covariance function as
\begin{eqnarray*}
\Gamma_e(s,t)&=&\cov\Big[\{\alpha Y\sum_{i=1}^\infty \frac 1{i^{2+\delta}}
\psi_i(s)\} ,\{\alpha Y\sum_{i=1}^\infty \frac 1{i^{2+\delta}}
\psi_i(t)\}\Big] \\
&=&\alpha^2\Big[\sum_{i=1}^\infty \frac 1{i^{2+\delta}}
\psi_i(s)\Big] \Big[\sum_{i=1}^\infty \frac 1{i^{2+\delta}}  \psi_i(t)\Big].
\end{eqnarray*}

\begin{Pro}\label{pro:my}
The following two optimization problems  \bse \argmax_\beta \frac
{\langle \Gamma_e \beta, \beta \rangle_{L_2}} {\langle \Gamma \beta,
\beta \rangle_{L_2}} \equiv \argmax_\beta \frac {\langle \Gamma_e
\beta, \beta \rangle_{L_2}} {\langle \Gamma_w \beta, \beta
\rangle_{L_2}}. \ese have the same solution $\beta$ given by \be
\beta(t)= c \sum_{i=1}^\infty \frac 1{i^{\delta}}\, \psi_i(t),
\label{eq:Grenander} \ee for any constant $c$.
\end{Pro}
\begin{proof}
From $\Gamma=\Gamma_e+\Gamma_w$, we have
\[
\frac {\langle \Gamma \beta, \beta \rangle_{L_2}}{\langle \Gamma_e \beta, \beta \rangle_{L_2}}=\frac{\langle \Gamma_e \beta, \beta \rangle_{L_2}+\langle \Gamma_w \beta, \beta \rangle_{L_2}}{\langle \Gamma_e \beta, \beta \rangle_{L_2}}
=1+\frac{\langle \Gamma_w \beta, \beta \rangle_{L_2}}{\langle \Gamma_e \beta, \beta \rangle_{L_2}}.
\]
Therefore, \bse \argmax_\beta \frac {\langle \Gamma_e \beta, \beta
\rangle_{L_2}} {\langle \Gamma \beta, \beta \rangle_{L_2}} \equiv
\argmax_\beta \frac {\langle \Gamma_e \beta, \beta \rangle_{L_2}}
{\langle \Gamma_w \beta, \beta \rangle_{L_2}}. \ese Let
$\beta=\sum_i b_i \psi_i$. Then \bse
\Gamma_e \beta &=& \alpha^2 \left(\sum_{i=1}^\infty \frac {b_i}{i^{2+\delta}}\right)\left( \sum_{i=1}^\infty \frac 1 {i^{2+\delta}}\psi_i\right),\\
\langle \Gamma_e \beta, \beta \rangle_{L_2}&=&\alpha^2 \left(\sum_{i=1}^\infty \frac {b_i} {i^{2+\delta}}\right)^2,\\
\Gamma_w \beta &=& \sum_{i=1}^\infty \frac {b_i} {i^{2}}\psi_i,\\
\langle \Gamma_w \beta, \beta \rangle_{L_2}&=&\sum_{i=1}^\infty
\frac {b_i^2} {i^{2}}. \ese Therefore, the optimization problem
becomes to maximize \bse \frac {\alpha^2 \left(\sum_{i=1}^\infty
{b_i}/ {i^{2+\delta}}\right)^2}{\sum_{i=1}^\infty {b_i^2}/ {i^{2}}}.
\ese From Cauchy-Schwarz inequality, \bse \left(\sum_{i=1}^\infty
\frac {b_i} {i^{2+\delta}}\right)^2\leq\left(\sum_{i=1}^\infty \frac
{b_i^2} {i^{2}}\right)\left(\sum_{i=1}^\infty \frac 1
{i^{2+2\delta}}\right). \ese The equality holds when $b_i \propto
1/i^\delta$, which means \bse \beta(t) \propto \sum_{i=1}^\infty
\frac 1{i^{\delta}}\, \psi_i(t)\ese is the maximum eigenfunction.
\end{proof}

The dimension reduction function $\beta(t)$ is obtained from solving
the eigenvalue problem $\Gamma_e\beta=\lambda \Gamma\beta$. The
corresponding optimal linear classification rule is via
\begin{equation}\label{lda_rule}
{\rm sign} \left(\langle \beta, X\rangle_{L_2}\right).
\end{equation}
This result can be linked to some prior study of linear discriminant
analysis of two Gaussian measures on a separable Hilbert space by
Grenander~\cite{grenander1950} and Rao and
Varadarajan~\cite{rao1963}. Let \bse m_y(t)\equiv
E\{X(t)|Y=y\}=\alpha y\sum_{i=1}^\infty \frac 1{i^{2+\delta}}\,
\psi_i(t) =\Gamma_w^{1/2} \left(\alpha y\sum_{i=1}^\infty
\frac{1}{i^{1+\delta}}\psi_i\right)(t) \in R(\Gamma_w^{1/2}). \ese
Note that $\beta$ given in~(\ref{eq:Grenander}) is not in $L^2(I)$,
but in $R(\Gamma_w^{-1/2})$, since
\[\|\Gamma_w^{1/2}\beta\|^2_{L_2}= \Big \|\sum_i
\frac{1}{i^{\delta+1}}\psi_i(t)
\Big\|^2_{L_2}=\sum_i\frac{1}{i^{2+2\delta}}<\infty.\] We also have
$\|\Gamma_w^{-1/2}m_y\|_{L_2}^2=\alpha^2\sum_{i=1}^\infty\frac{1}{\delta^{2+2\delta}}<\infty$,
i.e., $m_y$ is in $R(\Gamma_w^{1/2})$. Furthermore, from
Proposition~\ref{pro:my} and its proof, we have $\Gamma_e \beta =
c_1m_y$, where $c_1=c\alpha y
   \sum_{i=1}^\infty 1/(i^{2+2\delta})$,
 $\Gamma_w \beta =  c_2m_y$, where $c_2=c(\alpha y)^{-1}$.
Therefore
\begin{eqnarray*}
&& \langle \Gamma^{1/2} \beta, \Gamma^{1/2} \beta\rangle_{L_2}
 = \langle \Gamma \beta, \beta \rangle_{L_2}\\
&=& \langle (\Gamma_w+\Gamma_e) \beta, \beta \rangle_{L_2}
  = \langle \Gamma_w \beta, \beta \rangle_{L_2}+ \langle \Gamma_e \beta, \beta \rangle_{L_2}\\
&=& \|\Gamma_w^{1/2} \beta\|_{L_2}^2 + c_1c_2\langle m_y, \Gamma_w^{-1} m_y \rangle_{L_2}\\
&=& \|\Gamma_w^{1/2} \beta\|_{L_2}^2 + c^2 \left(\sum_{i=1}^\infty
  \frac{1}{i^{2+2\delta}}\right)\| \Gamma_w^{-1/2} m_y\|^2_{L_2}<\infty.
\end{eqnarray*}
Hence, $\beta \in R(\Gamma^{-1/2})$.

This is an example that $X\in L_2(I)$, $\beta\in R(\Gamma^{-1/2})$,
but $\beta\notin L_2(I)$, and the classification rule ${\rm sign}
\left(\langle \beta, X\rangle_{L_2}\right)$ is well-defined. This
indicates that,} to solve for a linear discriminant analysis problem
in $L_2(I)$, we cannot restrict $\beta$ to $L_2(I)$. We are obliged
to enlarge the domain of $\beta$ to $R(\Gamma^{-1/2})$. On the other
hand, requiring $\beta\in R(\Gamma^{-1/2})$ is indeed sufficient for
the purpose of linear discriminant analysis given
in~(\ref{lda_rule}) for classifying the observations into two
groups.
\end{Ex}

\begin{Ex}[Categorical response]
{\rm The feature revealed in Example 1 is not unique for binary
response variable $Y$. When the response variable $Y$ is
categorical, similar phenomenon can be observed. For example,
consider the case, where the response variable $Y$ is categorical
with possible values $y_1, \dots, y_k$. We normalize the $y$ values
so that $Y$ has mean zero and variance 1. Let
\begin{equation}
X_y (t) = \alpha y\sum_{i=1}^\infty \frac 1{i^{2+\delta}}  \psi_i(t)
 + \sum_{i=1}^\infty \frac 1i \, Z_i\, \psi_i(t),~~ t\in I.
\end{equation}
We can easily verify that the within-group covariance function is
\bse \Gamma_w(s,t)\equiv E\left[\cov\{X(s), X(t)|Y\}\right]
=\sum_{i=1}^\infty \frac 1 {i^2} \psi_i(s) \psi_i(t), \ese and the
between-group covariance function is \bse \Gamma_e(s,t)&\equiv&
\cov\left\{\alpha Y\sum_{i=1}^\infty \frac 1{i^{2+\delta}}  \psi_i(s), \alpha Y\sum_{i=1}^\infty \frac 1{i^{2+\delta}}  \psi_i(t)\right\}\\
&=&\alpha^2\Big\{\sum_{i=1}^\infty \frac 1{i^{2+\delta}}\,
\psi_i(s)\Big\} \Big\{\sum_{i=1}^\infty \frac 1{i^{2+\delta}}\,
\psi_i(t)\Big\}, \ese
 Let $\Gamma=\Gamma_w+\Gamma_e$. Note that the
forms of $\Gamma_w(s,t)$ and $\Gamma_e(s,t)$  here are exactly the
same as those in Example~1. Thus, when we perform the functional
sliced inverse regression by solving for the first eigenfunction,
\bse \beta_1=\argmax_v \frac {\langle \Gamma_e \beta, \beta
\rangle_{L_2}} {\langle \Gamma_w \beta, \beta \rangle_{L_2}}, \ese
we have exactly the same analysis as that in Example~1. It then
leads to the same conclusion. That is, we are obliged to enlarge the
domain of $\beta$ to $R(\Gamma^{-1/2})$. On the other hand,
requiring $\beta\in R(\Gamma^{-1/2})$ is also sufficient for our
purpose of classifying the observations into $k$ groups.}
\end{Ex}

\begin{Ex}[Continuous response]
{\rm Finally we provide an example with continuous response variable
$Y$. Let $Y$ have mean zero and variance 1, and let
\begin{equation}
X_y (t) = \alpha y\sum_{i=1}^\infty \frac 1{i^{2+\delta}}  \psi_i(t)
 + \sum_{i=1}^\infty \frac 1i \, Z_i\, \psi_i(t),~~ t\in I.
\end{equation}
We can easily verify that
the within-group covariance function is \bse \Gamma_w(s,t)\equiv
E\cov\{X(s), X(t)|Y\} =\sum_{i=1}^\infty \frac 1 {i^2} \psi_i(s)
\psi_i(t). \ese The between-group covariance function is \bse
\Gamma_e(s,t)&\equiv&
\cov\left\{\alpha Y\sum_{i=1}^\infty \frac 1{i^{2+\delta}}  \psi_i(s), \alpha Y\sum_{i=1}^\infty \frac 1{i^{2+\delta}}  \psi_i(t)\right\}\\
&=&\alpha^2 \Big\{\sum_{i=1}^\infty \frac 1{i^{2+\delta}}\,
\psi_i(s)\Big\} \Big\{\sum_{i=1}^\infty \frac 1{i^{2+\delta}}\,
\psi_i(t) \Big\}. \ese Let $\Gamma=\Gamma_w+\Gamma_e$. Now the same
analysis as that in Examples~1 and~2 leads to the conclusion that,
regardless of how many slices one decides to use, $\beta$ is in
$R(\Gamma^{-1/2})$.}
\end{Ex}

\section{Enlarged dimension reduction space and main results}
\label{sec:main}

In this section, we present our main results. First, we establish in
Theorem~\ref{th:one} an interesting link between covariance
operators on $L_2(I)$ and on $\rkhsg$. Next, we extend the
functional dimension reduction to a relaxed model with enlarged
space given in~(\ref{eq:fsir_model}). The reproducing kernel Hilbert
space $\rkhsg$, induced from the covariance operator $\Gamma$,
defines a proper range space for the sliced mean (see
Proposition~\ref{pro:5} and Theorem~\ref{th:mean}(a) below). It also
plays the parallel role as the span of $X$ in finite dimension (see
Proposition~\ref{pro:5}). Note that, $\rkhsg$ is equipped with an
inner product $\langle\cdot,\cdot\rangle_{\rkhsg}$. Interestingly,
this inner product  refers to the standardization (see
equation~(\ref{rkhs_norm}) above and equation~(\ref{rkhs_id}) below)
similar to the Mahalanobis distance and the standardization by the
covariance matrix in finite vector case. We also study the linear
design condition under the relaxed model in
Proposition~\ref{prop:7}.

\subsection{Bounded operators on $L_2(I)$ and on $\rkhsg$}

\begin{Th}\label{th:one}
Assume $\Gamma$ and $\Gamma_e$ are continuous, and respectively
strict positive definite and positive semi-definite. Then,
$\Gamma^{-1/2}\Gamma_{e}\,\Gamma^{-1/2}$ is a well-defined bounded
linear operator on $L_2(I)$ if and only if $\Gamma_{e}$ is a
well-defined bounded linear operator on $\rkhsg$.
\end{Th}
\begin{proof}
Let $h=\sum_i c_{i}\phi_i$. Then, \be \label{eq:norm_eq}
\|\Gamma^{-1/2}h\|_{L_2}^2 = \|\Gamma^{-1/2}\sum_i
c_{i}\phi_i(\cdot)\|_{L_2}^2 =\sum_{i}c_{i}^2/\xi_i=
\|h\|_{\rkhsg}^2. \ee That is, \be\label{eq:equiv} \Gamma^{-1/2}
h\in L_2(I) \Leftrightarrow h\in\rkhsg. \ee For any $g\in L_2(I)$,
there exists $h=\Gamma^{1/2}g\in \rkhsg$. Then
\[
\|\Gamma^{-1/2} \Gamma_e \Gamma^{-1/2}g\|_{L_2}=\|\Gamma^{-1/2} \Gamma_e h\|_{L_2}=\|\Gamma_e h\|_{\rkhsg}.
\]
Together with (\ref{eq:norm_eq}), we have
\[
\frac {\|\Gamma^{-1/2} \Gamma_e \Gamma^{-1/2}g\|_{L_2}}{\|g\|_{L_2}^2}=\frac{\|\Gamma_e h\|_{\rkhsg}}{\|h\|_{\rkhsg}^2},
\]
which yields the statement of the theorem.
\end{proof}

\begin{Rem}\label{rem:1}
From (\ref{eq:norm_eq}), $\Gamma^{-1/2}$ is bounded when it is
defined as a linear operator from $\rkhsg$ to $L_2(I)$. Here
boundedness is referred to its induced operator norm,
$\sup_{f\in\rkhsg} \|\Gamma^{-1/2} f\|_{L_2}^2/\|f\|_{\rkhsg}^2
<\infty$. However, when it operates on $f\in L_2(I)$,
$\Gamma^{-1/2}f$ may not belong to $L_2(I)$. For example,
$\Gamma^{-1/2}\phi_i(t) =\xi_i^{-1/2}\phi_i(t)$, hence
${\|\Gamma^{-1/2}\phi_i\|_{L_2}^2}/{\|\phi_i\|_{L_2}^2}=\xi_i^{-1}\to\infty$,
as $i\to\infty$. When combined with the additional covariance
operator $\Gamma_e$, Theorem~\ref{th:one} ensures the resulting
operator $\Gamma^{-1/2}\Gamma_e\Gamma^{-1/2}$ is a bounded linear
operator on $L_2(I)$, i.e.,
$\Gamma^{-1/2}\Gamma_e\Gamma^{-1/2}:L_2(I)\mapsto L_2(I)$ is  a
well-defined bounded operator. Note that $L_2(I)$ is a much larger
space than $\rkhsg$. Thus, the new operator composed of the three
operators can be well-defined in a larger domain than the original
operator $\Gamma^{-1/2}$ can.
\end{Rem}

\subsection{Relaxed model and extended estimation}

We are now in a position to revisit the functional dimension
reduction problem studied in~\cite{jiangetal2014}, describe the
problem more rigorously and extend it. Let $X(t)$, $t\in I$, be a
stochastic process satisfying $E\int X^2(t)dt <\infty$. Denote its
covariance function and spectrum by
 \bse \Gamma(s,t)\equiv \cov\{X(s),X(t)\}=\sum_{i=1}^\infty
\xi_i\phi_i(s)\phi_i(t). \ese Then, $X$ can be expressed by an
expansion as \bse X(s) =\sum_{i=1}^\infty A_i\phi_i(s), \ese where
$A_i$'s are independent random variables with mean zero and
variances $\xi_i$'s. Below we give a Proposition, which ensures that
we can exchange the order of double integrals.
\begin{Pro}\label{prop:exchange}
\[
E \langle X, \phi_i \rangle_{L_2} = \langle E(X), \phi_i \rangle_{L_2}.
\]
\end{Pro}
\begin{proof}
From Cauchy-Schwarz inequality, we have
\[
E \int |X(s) \phi_i(s) |ds \leq E\left[\left(\int X^2(s)ds \right)^{1/2} \left(\int \phi_i^2(s)ds \right)^{1/2}\right] =E\left(\int X^2(s)ds \right)^{1/2}.
\]
From Jensen's inequality,
\[
E\left(\int X^2(s)ds \right)^{1/2} \leq  \left(E\int X^2(s)ds \right)^{1/2} < \infty.
\]
Thus, with $E \int |X(s) \phi_i(s) |ds < \infty$, we can apply
Fubini's Theorem and get
\[
E \int X(s) \phi_i(s) ds = \int E\left[X(s)\right] \phi_i(s) ds.
\]
\end{proof}

Our proposed model is
\begin{equation}\label{eq:fsir_model}
Y=f\left(\langle \beta_1,X\rangle_{L_2},\dots, \langle
\beta_d,X\rangle_{L_2},\epsilon\right),~~ {\rm where}~~ \beta(\cdot)\in
R(\Gamma^{-1/2}).
\end{equation}
Note that a critical difference of our formulation here from that in
\cite{jiangetal2014} is that, we only require $\beta$ to be in
$R(\Gamma^{-1/2})$, which is larger than $L_2(I)$. This extension
allows more flexibility in the dimension reduction functions.

\begin{Pro}\label{prop:beta_space}
For $\beta\in R(\Gamma^{-1/2})$, $\langle \beta,X\rangle_{L_2}$ is
well-defined almost surely.
\end{Pro}

\begin{proof}
Let $\delta \equiv \Gamma^{1/2}\beta\in L_2(I)$ and $\delta_i
=\langle \delta,\phi_i\rangle_{L_2}$. We have \bse
&&E\left(\langle \beta, X\rangle_{L_2}\right)^2\\
&=&E\left( \sum_i \langle \Gamma^{-1/2}\delta, \phi_i\rangle_{L_2}
\cdot
\langle X,\phi_i \rangle_{L_2}\right)^2\\
&=&\sum_i\left(\xi_i^{-1/2} \delta_i \right)^2E A_i^2\\
&=& \sum_i \xi_i^{-1} \delta_i^2 \xi_i
  = \sum_i\delta_i^2
  = \|\delta\|_{L_2}^2<\infty,
\ese which implies that $|\langle \beta, X\rangle_{L_2}|<\infty$
a.s.
\end{proof}
\begin{Rem}\label{rem:Xspace}
Proposition \ref{prop:beta_space} reveals an interesting
result regarding the space where $X$ belongs to.
The finite second moment condition is commonly used
in statistical analysis.
In the finite dimensional case, a random
vector with finite second moment can have arbitrary variation
for each component of the random vector, hence the random vector
can take values in the entire space.
However, this is not the case
in the infinite dimensional functional space.
To ensure finite integrated variance,
a random function cannot have arbitrary variation
along each dimension. In fact, the variations along all dimensions,
except a finite set of dimensions, have
to degenerate sufficiently fast to guarantee finite total variant.
In fact, the set of dimensions in
which almost all variation accumulate is fixed for a single random
function. As
a consequence, the random function cannot
take values everywhere in $L_2(I)$. This is why
the resulting space of the random function $X$ is in fact a much smaller
subspace of $L_2(I)$. A feature of this subspace is that it ensures
finite inner-product with elements in $R(\Gamma^{-1/2})$, where
$\Gamma$ is the covariance function of $X$. We define this space as
\bse R(\Gamma^{1/2})^+ \equiv \{f: \langle
f,\beta\rangle_{L_2}<\infty, a.s. \quad \forall
 \beta \in R(\Gamma^{-1/2})\}.
\ese
Obviously $R(\Gamma^{1/2})\subset R(\Gamma^{1/2})^+\subset
L_2(I)$. We will encounter this space again when we present an equivalent
linearity condition later in Section~\ref{sec:linear}.
Note that although a single random function $X$ belongs to a much
smaller space $R(\Gamma^{1/2})^+$, the (uncountable) union of all such spaces of all random
functions  is the entire $L_2(I)$.
\end{Rem}

\begin{Rem}\label{rem:normshift}
 For any $f\in R(\Gamma^{1/2})$,
 Proposition~\ref{prop:beta_space} ensures that
the quantity $\langle \Gamma^{-1}f, X\rangle_{L_2}$ is well-defined
a.s. It is easy to verify  the  identity
\begin{equation}\label{rkhs_id}
\langle \Gamma^{-1/2}f, \Gamma^{-1/2}X\rangle_{L_2}=\langle \Gamma^{-1}f, X\rangle_{L_2}=\langle f, X\rangle_{\rkhsg}.
\end{equation}
In the classical SIR, the main problem can be viewed as solving the
eigenvalue problem of $\Gamma_e$ in the space scaled by
$\Gamma^{-1/2}$. Now in Functional Sliced Inverse Regression (FSIR),
(\ref{rkhs_id}) indicates that $\Gamma^{-1/2}$ can be again viewed
as the scaled operator from $L_2(I)$ to $\rkhsg$.
\end{Rem}

In fact, the relaxed model leads to more flexible requirements on
subsequent operators needed in the estimation procedure, which in
turn leads to less stringent conditions on quantities such as mean
covariates conditional on the response, etc. For example, in the
FSIR approach, we would search for $\beta$ from the functional
eigenvalue problem \be\label{eq:eig1} \Gamma_e\beta=\lambda
\Gamma\beta, \ee where $\Gamma(s,t)\equiv\cov\{X(s),X(t)\}$ as
before, $m_Y(s)\equiv E\{X(s)\mid Y\}$ and
\[\Gamma_e(s,t)\equiv \cov[E\{X(s)\mid Y\},E\{X(t)\mid Y\}] =\cov\{m_Y(s),m_Y(t)\}.\]
Letting $\eta=\Gamma^{1/2}\beta\in L_2(I)$, rewriting
(\ref{eq:eig1}) as \bse \Gamma^{-1/2}\Gamma_e\Gamma^{-1/2}\eta
=\Gamma^{-1/2}\Gamma_e\Gamma^{-1/2}(\Gamma^{1/2}\beta)
=\lambda(\Gamma^{1/2}\beta) =\lambda\eta, \ese we would naturally
require $\Gamma^{-1/2}\Gamma_e\Gamma^{-1/2}$ to be a well-defined
operator from $L_2(I)$ to $L_2(I)$. However,
$\Gamma^{-1/2}\Gamma_e\Gamma^{-1/2}$ is restricted to operate on
$R(\Gamma^{1/2})$ in \cite{jiangetal2014}. This restriction, comes
naturally from their condition that $\beta\in L_2(I)$, leads to a
conclusion that the slice mean can only be in a restricted space
$R(\Gamma)$ (Theorem~\ref{th:mean}(b) below) instead of in the space
$R(\Gamma^{1/2})$. In Theorem~\ref{th:mean}(a), we show that our
relaxation on the domain of $\Gamma^{-1/2}\Gamma_e\Gamma^{-1/2}$
leads to a more flexible condition on the conditional mean functions
$m_Y(s)$.

Here, we first state a useful result in Proposition \ref{pro:5}.

\begin{Pro}\label{pro:5}
$R(\Gamma^{1/2})\equiv \rkhsg$.
\end{Pro}
\begin{proof} A function $g\in R(\Gamma^{1/2})$ is equivalent to
$g=\Gamma^{1/2}h$ and $h\in L_2(I)$. Now \bse
\|g\|^2_{\rkhsg}&=&\|\Gamma^{1/2}h\|_{\rkhsg}^2\\
&=&\|\sum_i \xi_i^{1/2}\phi_i(s)\langle\phi_i(t),h(t)\rangle\|_{\rkhsg}^2\\
&=&\|\sum_i \{\xi_i^{1/2}\langle\phi_i(t), h(t)\rangle\}^2/\xi_i\\
&=&\|h\|_{L_2}^2. \ese Thus $\|g\|^2_{\rkhsg}<\infty$ is equivalent
to $\|h\|_{L_2}^2<\infty$, hence $g\in R(\Gamma^{1/2})$ is
equivalent to $g\in\rkhsg$.
\end{proof}

\begin{Rem}
Proposition~\ref{pro:5} implies that, if $m_y\in\rkhsg$, then
$m_y\in R(\Gamma^{1/2})$, and thus $\Gamma^{-1}m_y$ is in
$R(\Gamma^{-1/2})$. In those examples in Section~\ref{sec:ex},  we
have shown that the relaxation from $\beta\in L_2(I)$ to $\beta\in
R(\Gamma^{-1/2})$ is crucial and that the condition $\beta\in
R(\Gamma^{-1/2})$ is sufficient for $\langle \beta,X\rangle_{L_2}$
being well-defined a.s. Furthermore, from the proof of
Proposition~\ref{prop:beta_space}, we have \bse |\langle \beta, m_y
\rangle_{L_2}| =|\langle \beta, E(X|Y=y) \rangle_{L_2}|=|E(\langle
\beta, X \rangle_{L_2}|Y=y)|\le \left[E(\langle \beta, X
\rangle_{L_2}^2|Y=y)\right]^{1/2}<\infty,\ese which means that $m_y
\in R(\Gamma^{1/2})$. Therefore, our relaxed condition on $\beta$ is
sufficient to include all possible $m_y$. In fact, it is also
necessary since for any proper subset $\Omega \subsetneq
R(\Gamma^{-1/2})$, there always exists some $m_y$ so that the
optimal $\beta=\Gamma^{-1}m_y \notin \Omega$.

\end{Rem}

\begin{Th}\label{th:mean} Let $Y$ take values in a discrete finite
  set, say $\{1,\dots,k\}$,  with equal probability.\\
(a) If $\Gamma^{-1/2}\Gamma_e\Gamma^{-1/2}$ is a bounded operator
from
$L_2(I)$ to $L_2(I)$, then $m_y\in R(\Gamma^{1/2})$.\\
(b) Alternatively, if $\Gamma^{-1/2}\Gamma_e\Gamma^{-1/2}$ is a
bounded operator from $R(\Gamma^{1/2})$ to $R(\Gamma^{1/2})$, then
$m_y\in R(\Gamma)$.
\end{Th}
\begin{proof}
(a) From Theorem~\ref{th:one}, if
$\Gamma^{-1/2}\Gamma_e\Gamma^{-1/2}$ is a bounded operator on
$L_2(I)$, then $\Gamma_e$ is a bounded operator on $\rkhsg$, which
means $h\equiv\Gamma_e\, g \in \rkhsg$ for any $g \in \rkhsg$. Thus,
\begin{eqnarray*}
h=\Gamma_e\, g = \frac1k\sum_{y=1}^k m_y \otimes m_y\trans \, g
= \frac1k \sum_{y=1}^k \langle m_y, g \rangle_{L_2} \, m_y.
\end{eqnarray*}
In order for the above function to be in $\rkhsg$ for arbitrary
$g\in\rkhsg$, $m_y$'s have to be in $\rkhsg$.

(b) For an arbitrary $g \in R(\Gamma^{1/2})$,
 $g_1\equiv
\Gamma^{-1/2}g$ is in $L_2(I)$. Since $h\equiv
\Gamma^{-1/2}\Gamma_e\Gamma^{-1/2} g$ is in $R(\Gamma^{1/2})$, we
have $\Gamma^{1/2}h$ is in $R(\Gamma)$. Furthermore, $\Gamma^{1/2}h$
can be expressed as
\[
\Gamma^{1/2}h =\Gamma_e \Gamma^{-1/2}\, g
= \frac 1k \sum_{y=1}^k m_y \otimes m_y\trans  g_1= \frac 1k \sum_{y=1}^k\langle m_y, g_1 \rangle_{L_2}\, m_y.
\]
In order for the above function to be in $R(\Gamma)$ for arbitrary
function $g_1 \in L_2(I)$, $m_y$'s have to be in $R(\Gamma)$.
\end{proof}

We now examine how the formulation will affect the estimation
procedure. Assume a discrete $Y$ for simplicity. The $j^{\rm th}$
slice mean function is given by
\[m_j(t)=E\{X(t)|Y=j\}= \sum_i E(A_i|Y=j)\phi_i(t).\]
Following Theorem \ref{th:mean} and Proposition \ref{pro:5}, $m_j$
is in RKHS $\rkhsg$. Assume in the $j^{\rm th}$ slice, we have
observations ${\cal D}_j\equiv \{X_i({\cal T}),Y_i\}_{i=1}^{n_j}$,
where $Y_i=j$ and we consider two types of ${\cal T}$. One is ${\cal
T}=I$, i.e., we observe the whole sample paths, and the other is
${\cal T}=\{t_k\}_{k=1}^{q}, q<\infty$, i.e., we observe $X_i(t)$ at
some common discrete time points.

\begin{Th} [Representer Theorem]\label{thm:rep} Given the $j^{\rm th}$ slice training sample ${\cal D}_j$
with ${\cal T}=\{t_k\}_{k=1}^{q}$, an arbitrary empirical risk
function ${\mathbb Q}:\real^2\mapsto\real$ and a scalar $C>0$,
consider the following minimization problem
\begin{equation}\label{eq:minimization}\argmin_{m\in\rkhsg}
\sum_{i=1}^{n_j}\sum_{t\in{\cal T}} {\mathbb Q}\left\{X_i(t),m(t)\right\}
+C\|m\|^2_{\rkhsg}.\end{equation}
Then the solution of the minimization problem exists and has the
representation  form
\begin{equation}\label{eq:slice_mean}
{\widehat m}_j(s) =\sum_{k=1}^{q} \Gamma(s,
t_k)\alpha_{jk},~~ t_k\in{\cal T},~\alpha_{jk}\in\real.
\end{equation}
\end{Th}

\begin{proof}
For any function $m(s)\in \rkhsg$, it can be expressed as
\[
m(s)=\sum_{k=1}^q  \Gamma(s,t_k)\alpha_k + \nu(s),
\]
where $\nu(s)$ is in $\rkhsg$ and orthogonal to every
$\Gamma(s,t_{k})$, $j=1,\dots,k$. By the reproducing property of
$\rkhsg$,
\[
m(t_{\ell})=\left\langle \Gamma(s,t_{\ell}), \sum_{k=1}^q  \Gamma(s,t_{k})\alpha_k
 + \nu(s) \right\rangle_{\rkhsg}
 = \sum_{k=1}^q \Gamma(t_\ell, t_{k})\, \alpha_k,~~ \ell=1,\dots,q,
\]
which does not involve $\nu(s)$. This implies that the empirical
risk function ${\mathbb Q}$ in (\ref{eq:minimization})  also does
not involve $\nu(s)$. Since
\[
 \Big\|\sum_{k=1}^q   \Gamma(s,t_{k})\alpha_k + \nu(s)\Big\|^2_{\rkhsg}=
  \Big\|\sum_{k=1}^q   \Gamma(s,t_{k})\alpha_k \Big\|^2_{\rkhsg}+ \Big\|\nu\Big\|^2_{\rkhsg}
\geq \Big\|\sum_{k=1}^q   \Gamma(s,t_{k})\alpha_k\Big \|^2_{\rkhsg},
\]
the regularization term in~(\ref{eq:minimization}) is minimized by
$\nu(s)=0$. Therefore, the minimizer takes the form ${\widehat
m}_j(s) =\sum_{k=1}^q  \Gamma(s, t_{k})\alpha_{k}$.
\end{proof}

From the proof of Theorem~\ref{thm:rep}, we can see that the role of
$C\|m\|^2_{\rkhsg}$ in (\ref{eq:minimization}) is to force  $\nu$ to
be zero and hence to guarantee a unique solution of the minimization
problem. If we set $C=0$, $\nu$ can be chosen freely as any function
orthogonal to $\Gamma(s,t_k)$'s and it will not affect the target
value in (\ref{eq:minimization}). This freedom occurs because
$\{\Gamma(s, t_k)\}_{k=1}^q$ do not span the whole $\rkhsg$. If such
freedom vanishes, for example this happens when the entire
$X_i(t)$'s are observed, then we no longer need to have the
$C\|m\|^2_{\rkhsg}$ term to induce uniqueness. An additional utility
of the ``penalty'' term $C\|m\|^2_{\rkhsg}$ is to regularize the
solution. It provides a balance between the best data fit evaluated
by the risk function and the variability of the solution.

\begin{Rem}
If we modify the minimization problem~(\ref{eq:minimization}) by
restricting the residing space of the slice mean to a smaller
subspace $R(\Gamma)$,
\begin{equation}
\argmin_{m\in R(\Gamma)}
\sum_{i=1}^{n_j}\sum_{t\in{\cal T}} {\mathbb
Q}\left\{X_i(t),m(t)\right\} +C\|m\|^2_{\rkhsg}.
\end{equation}
then the representation form~(\ref{eq:slice_mean}) might not be
valid anymore.
\end{Rem}

\begin{Rem}
When the observations are the entire paths, i.e., ${\cal T}=I$, we
can choose $C= 0$ and modify the
minimization~(\ref{eq:minimization}) to
\begin{equation}\label{eq:minimization2}
\argmin_{m\in\rkhsg}
\sum_{i=1}^{n_j} {\mathbb Q}(X_i,m),
\end{equation}
where now $\mathbb Q$ is a bivariate risk functional. A typical
bivariate risk functional is the quadratic one, i.e., ${\mathbb
Q}(f,g) = \langle \Lambda (f-g), f-g\rangle_{L_2}$, where $\Lambda$
is a symmetric strictly positive definite linear integral operator
with $\zeta_\ell$ and $\psi_\ell$ ($\ell=1, \dots, \infty$) as its
eigenvalues and eigenfunctions. In this case, ${\mathbb
Q}(f,g)=\sum_{\ell=1}^\infty {\zeta_\ell} (f_\ell-g_\ell)^2$ for
$f=\sum_{\ell=1}^\infty f_\ell \psi_\ell$ and
$g=\sum_{\ell=1}^\infty g_\ell \psi_\ell$. Write $X_i$ as
$X_i=\sum_{\ell=1}^\infty x_{i\ell} \psi_\ell$ and
$m=\sum_{\ell=1}^\infty m_\ell \psi_\ell$. Then
\begin{eqnarray*}
\sum_{i=1}^{n_j}{\mathbb Q}\left(X_i,m \right)
=\sum_{i=1}^{n_j} \sum_{\ell=1}^\infty {\zeta_\ell} (x_{i\ell}-m_\ell)^2
  =  \sum_{\ell=1}^\infty {\zeta_\ell} \sum_{i=1}^{n_j} (x_{i\ell}-m_\ell)^2.
\end{eqnarray*}
The above term is minimized when $m_\ell=\sum_{i=1}^{n_j}
x_{i\ell}/n_j$ for all $\ell$. That is, the mean path is the
minimizer of~(\ref{eq:minimization2}).
\end{Rem}

\begin{Rem}
With a given covariance estimator, the slice means can be expressed
as a linear combination of covariance functions at training data
points, as presented in~(\ref{eq:slice_mean}). When the covariance
estimator is given, the estimation of slice means becomes less
challenging. The most difficult part of estimation in FSIR is the
estimation of covariance operator. High-dimensional covariance
estimation is a difficult problem, and the functional case is even
more challenging. Our aim here is to set up a right framework for
the functional inverse regression in an enlarged space. Therefore,
we do not further discuss the estimation of the covariance operator.
\end{Rem}

\subsection{Linearity condition re-expressed}\label{sec:linear}

Recall that SIR requires a linearity condition, which, in the
functional dimension reduction framework, is written as the
following: For any $b\in R(\Gamma^{-1/2})$ there exist
$a_0,a_1,\dots,a_k\in\real$ such that
\begin{equation}\label{eq:ldc_gfsir0}
E\left(\langle b, X\rangle_{L_2} |\langle \beta_1,X\rangle_{L_2},\dots,
\langle \beta_k,X\rangle_{L_2}\right)
= a_0+ \sum_{j=1}^k a_j \langle \beta_j,X\rangle_{L_2},
\end{equation}
where $\beta_1, \dots, \beta_k\in R(\Gamma^{-1/2})$. Below we give a
more direct linearity condition statement, which is equivalent to
the one given by~(\ref{eq:ldc_gfsir0}). \bse E\left(X(s)|\langle
\beta_1,X\rangle_{L_2},\dots, \langle \beta_k,X\rangle_{L_2}\right)~
\mbox{is linear in } \langle \beta_1,X\rangle_{L_2},\dots, \langle
\beta_k,X\rangle_{L_2}, ~ \forall s\in I, \ese where $\beta_1,
\dots, \beta_k\in R(\Gamma^{-1/2})$. That is, there exist
$a_j(\cdot)'s\in R(\Gamma^{1/2})^+$ such that \be
\label{eq:ldc_gfsir} E\left(X(s)|\langle
\beta_1,X\rangle_{L_2},\dots, \langle
  \beta_k,X\rangle_{L_2}\right)
=a_0(s)+\sum_{j=1}^k a_j(s)\langle \beta_j,X\rangle_{L_2}, ~ \forall
s\in I. \ee

\begin{Rem}
 It is easy to check that $R(\Gamma^{1/2})^+ \subset L_2(I)$.
Since the functions $a_j$'s in~(\ref{eq:ldc_gfsir}) should belong to
the same space where $X$ resides, they belong to $R(\Gamma^{1/2})^+$
based on Proposition~\ref{prop:beta_space}, which  is smaller than
$L_2(I)$. This fact about $a_j$'s is masked when
(\ref{eq:ldc_gfsir0}) is used to describe the linearity condition.
However, we can see that this condition on $a_j$'s is indeed
necessary and sufficient from the following proof of
Proposition~\ref{prop:7}.
\end{Rem}

\begin{Pro}\label{prop:7}
The two versions of functional linearity condition given in
(\ref{eq:ldc_gfsir0}) and (\ref{eq:ldc_gfsir}) are equivalent.
\end{Pro}

\begin{proof}
Assume (\ref{eq:ldc_gfsir0}) holds. Consider the evaluation
functional ${\cal F}_s (X)=X(s)$.  Let $b_s(\cdot)\equiv
\sum_i\phi_i(s)\phi_i(\cdot)$. Since $\|\Gamma^{1/2}
b_s(\cdot)\|_{L_2}^2=\sum_i\xi_i\phi_i^2(s)=\Gamma(s,s)<\infty$ for
any $s$, we have $b_s(\cdot)\in R(\Gamma^{-1/2})$ for any $s$.
Obviously
\begin{equation}\label{Riesz_rep}
\langle b_s, X\rangle_{L_2}={\cal F}_s (X)=X(s).
\end{equation}
Thus, \be E\left\{X(s) |\langle \beta_1,X\rangle_{L_2},\dots,
\langle \beta_k,X\rangle_{L_2}\right\} &=&E\left(\langle b_s,
X\rangle_{L_2} |\langle \beta_1,X\rangle_{L_2},\dots,
\langle \beta_k,X\rangle_{L_2}\right)\nonumber\\
&=&  a_0(s) + \sum_{j=1}^k a_j(s) \langle \beta_j,X\rangle_{L_2}.
\label{eq:ldc_gfsir1} \ee By Proposition~\ref{prop:beta_space} we
have $X \in R(\Gamma^{1/2})^+$. It is then easy to see from the
identities~(\ref{eq:ldc_gfsir1}) that $a_j$'s are
in~$R(\Gamma^{1/2})^+$. Hence (\ref{eq:ldc_gfsir}) holds.

On the other hand, assume (\ref{eq:ldc_gfsir}) holds, \bse
E\left(X(s) |\langle \beta_1,X\rangle_{L_2},\dots, \langle
\beta_k,X\rangle_{L_2}\right) =a_0(s)+ \sum_{j=1}^k a_j(s) \langle
\beta_j,X\rangle_{L_2}. \ese Now for any $b(s)\in R(\Gamma^{-1/2})$,
take inner product with the above two sides, we obtain \bse
E\left(\langle b, X\rangle_{L_2} |\langle
\beta_1,X\rangle_{L_2},\dots, \langle \beta_k,X\rangle_{L_2}\right)
=\langle b, a_0\rangle_{L_2} + \sum_{j=1}^k\langle b, a_j
\rangle_{L_2} \langle \beta_j,X\rangle_{L_2}.\ese Since $b \in
R(\Gamma^{-1/2})$ and $a_j \in R(\Gamma^{1/2})^+$, $\langle b, a_j
\rangle_{L_2}<\infty$.  Therefore, (\ref{eq:ldc_gfsir0}) holds.
\end{proof}

\begin{Rem}\label{rem:linearCondition}
We provide a neat expression  (\ref{eq:ldc_gfsir}) for the linearity
condition.
 In the classical SIR, the corresponding condition of
(\ref{eq:ldc_gfsir0}) is:
 For any ${\bm b}\in \real^p$, there exist
$a_0,a_1,\dots,a_k\in\real$ such that
\[
E\left({\bm b}\trans  \X|\bb_1\trans  \X,\dots, \bb_k\trans  \X\right)
=a_0+\sum_{j=1}^k a_j \bb_j\trans  \X. \] The corresponding condition of
(\ref{eq:ldc_gfsir}) is: There exist ${\bm a}_j$'s in $\real^p$ such
that
\[
E\left(\X|\bb_1\trans  \X,\dots, \bb_k\trans  \X\right) ={\bm
a}_0+\sum_{j=1}^k {\bm a}_j \bb_j\trans  \X. \] They are equivalent
by similar arguments above. Interestingly, such equivalence
description of the functional linearity condition seems only
possible when we allow $\beta\in R(\Gamma^{-1/2})$. In the original
framework of \cite{jiangetal2014}, where $\beta$ is required to be
in $L_2(I)$, we are unable to obtain such equivalence description,
as the representation function $b_s(\cdot)$ in~(\ref{Riesz_rep}) for
the evaluation functional ${\cal F}_s$ is in $R(\Gamma^{-1/2})$ but
not in $L_2(I)$.
\end{Rem}

\section{Conclusion}\label{se:conclude}

We have described an extension of the dimension reduction models to
the functional data framework.  Our extension is based on careful
and rigorous considerations in operator theory and functional
analysis. We mainly focused on generalizing concepts in the
classical dimension reduction problems into the new framework and on
enlarging the functional space of the reduction function $\beta$. We
found some interesting examples where such increased flexibility is
indeed needed, and we discovered an equivalent expression of the
popular linearity condition. While our analysis is based on FSIR, we
believe similar analysis can be applied to other functional inverse
regression based methods. It will be interesting to study how other
methods in the classical dimension reduction models can be properly
extended to the functional data framework.

\bibliographystyle{apalike}
\bibliography{ref}
\end{document}